\documentclass[12pt]{article}
\usepackage{amsfonts}
\usepackage{amsmath,color,amsthm,amssymb,cite,bbm,mathrsfs}
\usepackage{epsfig}
\usepackage{graphicx}
\usepackage{latexsym}

\usepackage{color}
\title{Logarithmic Representability of Integers as $k$-Sums}
\author {Anant P.~Godbole\\
Department of Mathematics and Statistics\\
East Tennessee State University\and Samuel Gutekunst\\
Department  of Mathematics\\
Harvey Mudd College\and Vince Lyzinski\\ Department of Applied Mathematics and Statistics\\The Johns Hopkins University\and Yan Zhuang\\
Department of Mathematics and Computer Science\\
Goucher College}
\begin{document}
\def\qed{\vbox{\hrule\hbox{\vrule\kern3pt\vbox{\kern6pt}\kern3pt\vrule}\hrule}}
\def\ms{\medskip}
\def\n{\noindent}
\def\ep{\varepsilon}
\def\G{\Gamma}
\def\lr{\left(}
\def\ls{\left[}
\def\rs{\right]}
\def\lf{\lfloor}
\def\rf{\rfloor}
\def\lg{{\rm lg}}
\def\lc{\left\{}
\def\rc{\right\}}
\def\rr{\right)}
\def\ph{\varphi}
\def\p{\mathbb P}
\def\nk{n \choose k}
\def\cA{{\cal A}}
\def\s{\cal S}
\def\e{\mathbb E}
\newcommand{\begp}{\begin{proposition}}
\newcommand{\enp}{\end{proposition}}
\newcommand{\begt}{\begin{thm}}
\newcommand{\ent}{\end{thm}}
\newcommand{\begl}{\begin{lem}}
\newcommand{\enl}{\end{lem}}
\newcommand{\begc}{\begin{corollary}}
\newcommand{\enc}{\end{corollary}}
\newcommand{\begcl}{\begin{claim}}
\newcommand{\encl}{\end{claim}}
\newcommand{\begr}{\begin{remark}}
\newcommand{\enr}{\end{remark}}
\newcommand{\begal}{\begin{algorithm}}
\newcommand{\enal}{\end{algorithm}}
\newcommand{\begd}{\begin{definition}}
\newcommand{\enf}{\end{definition}}
\newcommand{\begx}{\begin{example}}
\newcommand{\enx}{\end{example}}
\newcommand{\bega}{\begin{array}}
\newcommand{\ena}{\end{array}}
\newcommand{\bsno}{\bigskip\noindent}
\newcommand{\msno}{\medskip\noindent}
\newcommand{\oM}{M}
\newcommand{\omni}{\omega(k,a)}
\newcommand{\yjn}{Y_{j,n}}
\newcommand\gd{\delta}
\newcommand\gD{\Delta}
\newcommand\gl{\lambda}
\newcommand\gL{\Lambda}
\newcommand{\ijn}{I_{j,n}}
\newcommand{\Bin}{B_{i,n}}
\newcommand{\yjkn}{Y_{j,k,n}}
\newcommand{\ijkn}{I_{j,k,n}}
\newcommand{\Bikn}{B_{i,k,n}}
\newcommand{\skan}{S_k(\alpha,n)}
\newcommand{\san}{S_2(\alpha,n)}
\hyphenation{Quick-sort}
\newcommand\urladdrx[1]{{\urladdr{\def~{{\tiny$\sim$}}#1}}}
\newtheorem{thm}{Theorem}[section]
\newtheorem{con}{Conjecture}[section]
\newtheorem{claim}[thm]{Claim}
\newtheorem{definition}[thm]{Definition}
\newtheorem{lem}[thm]{Lemma}
\newtheorem{cor}[thm]{Corollary}
\newtheorem{remark}[thm]{Remark}
\newtheorem{prp}[thm]{Proposition}
\newtheorem{ex}[thm]{Example}
\newtheorem{eq}[thm]{equation}
\newtheorem{que}{Problem}[section]
\newtheorem{ques}[thm]{Question}
\providecommand{\floor}[1]{\left\lfloor#1\right\rfloor}
\maketitle
\begin{abstract}  A set $\cA=\cA_{k,n}\subset[n]\cup\{0\}$ is said to be an additive $k$-basis if each element in $\{0,1,\ldots,kn\}$ can be written as a $k$-sum of elements of $\cA$ in {\it at least} one way.  Seeking multiple representations as $k$-sums,  and given any function $\phi(n)\to\infty$, we say that $\cA$ is said to be a  \emph{truncated $\phi(n)$-representative $k$-basis} for $[n]$ if 
for each $j\in[\alpha n, (k-\alpha)n]$ the number of ways that $j$ can be represented as a $k$-sum of elements of $\cA_{k,n}$ is $\Theta(\phi(n))$.  In this  paper, we follow tradition and focus on the case $\phi(n)=\log n$, and show that  a randomly selected set in an appropriate probability space is a truncated log-representative basis with probability that tends to one as $n\to\infty$.  This result is a finite version of a result proved by Erd\H os \cite{erdos} and extended by Erd\H os and Tetali \cite{ET}.

\end{abstract}
\section{Introduction}  In 1956 Erd\H os \cite{erdos} answered a question posed in 1932 by Sidon by proving that there exists an infinite sequence of natural numbers $\mathcal S$ and constants $c_1$ and $c_2$ such that for large $n$, \begin{equation}c_1\log n\leq r_2(n)\leq c_2\log n,\end{equation}
where, for $k\ge 2$, $r_k(n)$ is the {\it number of ways} of representing the integer $n$ as the sum of $k$ distinct elements from $\mathcal S$, a so-called {\it log-representative basis} of order $k$.  The result was generalized in the 1990 work of Erd\H os and Tetali \cite{ET} which established that there exists an infinite sequence $\mathcal S$ for which (1) was true for each fixed $k\geq 2$, i.e., for each large $n$, 
\begin{equation} r_k(n)=\Theta(\log n). 
\end{equation}
To achieve this result, Erd\H os and Tetali constructed a random sequence $\mathcal S$ of natural numbers by including $z$ in $\mathcal S$ with probability
$$p(z)=\begin{cases} C\frac{(\log{z})^{1/k}}{z^{(k-1)/k}}, & \mbox{if } z>z_0 \\ 0 & \mbox{otherwise } \end{cases}$$
where $C$ is a determined constant and $z_0$ is the smallest constant such that $p(z_0)\leq 1/2$.  They then showed that this random sequence is almost surely (a.s.)~a log-representative basis of order $k$, with (2) holding a.s.~for large $n$.  We note here that a.s.~in this context means ``with probability one," i.e., in the sense used in measure theory.

For a natural finite variant of the above problem, we define: 
\begd \emph{With $[n]:=\{1,2,\ldots,n\}$, a set $\cA_{k,n}\subseteq[n]\cup\{0\}$ is said to be a \emph{log-representative $k$-basis for B} (or simply a representative $k$-basis for $B$) if each $j\in B\subset[kn]\cup\{0\}$ can be represented as a $k$-sum of distinct elements of $\cA_{k,n}$ in $\Theta(\log n)$ ways. }\enf  
\begr
\emph{To see how this is the natural finite variant of the problem tackled in Erd\H os and Tetali in \cite{ET}, note that they teased out asymptotics for the emergence of log-representative bases for $B=[N_0,\infty)$ for some suitable $N_0$.  They showed that in some probability space, almost all infinite sequences $\mathcal S$ satisfy (2).  It is natural to then ask, for finite $\cA_{k,n},$ how small $\cA_{k,n}$ can be while still being a representative $k$-basis for a suitable $B$.}
\enr
\begr
\emph{Note that a more general definition might ask that the number of representations of $j$ equal $\Theta(\phi(n))$ for some $\phi(n)\to\infty$, but we stick close  to tradition and just deal with the case $\phi(n)=\log n$; it is interesting to note, though, that Sidon's original question asked about whether it was possible to find a representative basis with $\phi(n)=o(n^\epsilon)$ for all $\epsilon>0$.}
\enr

We will use a probability model in which each integer in $[n]\cup\{0\}$ is chosen to be in $\cA_{k,n}$ with equal (and low) probability $p=p_n$.  Since, e.g., the only way to represent 1 as a 2-sum of elements of $[n]\cup\{0\}$  is as 1+0, it would be impossible for the random ensemble to form an representative  basis unless we choose a target sumset, $B$, smaller than $[kn]\cup\{0\}$; this motivates the next definition -- which was the one adopted in \cite{AV} even when we had $\phi(n)=1$ for each $n$.

\begd  \emph{Let $\cA_{k,n}$ be a  subset of $[n]\cup\{0\}.$   Then $\cA_{k,n}$ is said to be a  \emph{truncated log-representative $k$-basis} for $[n]$ if 
each $j\in[\alpha n, (k-\alpha)n]$ can be represented as a $k$-sum of distinct elements of $\cA_{k,n}$ in $\Theta(\log n)$ ways. }
\enf

In \cite{AV}, the authors used Poisson approximation (see \cite{bhj} for background) and the Janson inequality \cite{janson} to derive a sharp threshold for which values of $p=p_n$ make the set $\cA_{k,n}$ almost never/almost surely  a truncated $k$-basis as $n\rightarrow\infty$, i.e.~every $j\in[\alpha n, (k-\alpha)n]$ could be represented \emph{at least once} as a $k$-sum of elements in $\cA_{k,n}$ with probability tending to 0 or 1 as $n\rightarrow\infty$.  Here the phrases ``almost never" and  ``almost surely" are used as is traditional in random methods.  The threshold function for $k$=2 (for example) is roughly at $p_n=A_\alpha{\sqrt{\log n/n}}$, with a third order correction term controlling the actual threshold, in contrast to the fact that the minimal size of a truncated 2-basis is of magnitude $O(\sqrt n)$ \cite{mrose}, \cite{hh}.   The corresponding questions of maximal Sidon families (i.e., ones for which each target integer is represented {\it at most once}), and zero-one thresholds for the emergence of the Sidon property feature a wider gap; see \cite{gjlr}.

The authors of \cite{AV}  did not derive the asymptotics for which $p$ determine whether $\cA_{k,n}$ is a truncated representative $k$-basis as $n\rightarrow\infty$, a question which we take up presently.  Our work is organized as follows:  We present results on truncated representative $2$-bases in Section 2, using some simple Chernoff bounds.
In Section 3, we consider similar questions for truncated representative $k$-additive bases, and we apply Talagrand's inequality \cite{as} to derive our desired results.

\begr
\emph{An alternate way of dealing with the boundary effects encountered in finite additive bases is to define \emph{modular representative $k$-bases.}   A set $\cA_{k,n}\subseteq[n-1]\cup\{0\}$ is said to be a modular representative $k$-basis for $[n]$ if the number of ways that each $j\in[n-1]\cup\{0\}$ can be written as a mod$(n)$ $k$-sum of elements of $\cA_{k,n}$ is $\Theta(\log n)$.  Definitive results on the emergence of modular additive bases have been proved in the papers of Yadin \cite{ya} using the method of Brun's sieve and in Sandor \cite {sandor} using Janson's correlation inequalities.  Neither tackled the representative basis question, as we do in the present work.  We believe, moreover, that the truncated basis is the more natural finite variant of the problem considered by Erd\H os and Tetali in \cite{ET}.  They were concerned with constructing a basis with $r_k(n)=\Theta(\log n)$ for all integers greater than a fixed but arbitrary $N_0$.  This allowed them great flexibility in choosing the threshold $N_0$ to be large enough to achieve the desired behavior.  One might ask how small $N_0$ can be while still maintaining an additive basis with $r_k(n)=\Theta(\log n)$, which is the natural analogue to the truncated basis question explored below.}
\enr

Throughout the rest of the paper, we suppress the descriptors  ``truncated" and ``log", referring simply to ``representative $k$-bases."

\section{2-Additive Representations}
Consider first the case where $k=2$.  Construct the random set $\mathcal{A}_{2,n}$ by choosing each integer in $[n]\cup\{0\}$ to be in $\cA_{2,n}$ independently with probability $p=p_n$.  Let $S_2(\alpha,n):=[\alpha n, (2-\alpha)n]$, and for each $j\in S_2(\alpha,n)$, let $Y_{j,n}$ be the number of ways that $j$ can be represented as a $2$-sum of distinct elements of $\cA_{2,n}$; the case where summands are allowed to be equal can be proved exactly as in what follows.  Let $I_{j,n}:=\mathbbm{1}\{Y_{j,n}\neq \Theta(\log n)\}$, so that $X_n:=\sum_{j\in S_2(\alpha,n)} I_{j,n}$ is the number of elements of $S_2(\alpha,n)$ that are not represented order $\log n$ times in the 2-sum set.  

For each $j\in [\alpha n, n],$ the maximum number of representations as 2-sums from $\cA_{2,n}$ is given by
$$\rho_{2,n}(j)=\rho_{2,n}(2n-j)=\left\lceil\frac{j}{2}\right\rceil.$$  Fixing $j$, for $i=1,...,\rho_{2,n}(j)$ let 
$$\Bin=\mathbbm{1}\{i\text{-th pair of integers in }[n]\cup\{0\}\text{ summing to {\it j} is present in }\cA_{2,n}\}$$  so that 
$\yjn=\sum_{i=1}^{\rho_{2,n}(j)}\Bin$.  Note that each integer in $[n]\cup\{0\}$ can be in at most one of the $\rho_{2,n}(j)$ pairs of integers summing to $j$, and so the associated $\Bin$'s are independent Bern$(p^2)$ random variables.  It follows that $\yjn$ has a  Bin($\rho_{2,n}(j),p^2$) distribution.  Two straightforward applications of Chernoff-type bounds then yield the following result:
\begt
\label{t:2rep}
Let $\alpha\in(0,1)$ be fixed, and let $\eta>0$ be arbitrarily small.  Create the random set $\cA_{2,n}$ by picking each integer in $[n]\cup\{0\}$ to be in $\cA_{2,n}$ with probability
$$p=p_n:=\sqrt{\frac{\left(\frac{2}{\alpha}+\eta\right)\log n}{n}}.$$ 
Then 
$$\lim_{n\rightarrow \infty}\p(X_n=0)=1,$$ so that w.h.p. $\cA_{2,n}$  is an asymptotic representative 2-basis as $n\rightarrow\infty$.
\ent
\begin{proof}  
First note that $\rho_{2,n}(j)$ is maximized by $j=n$, so that for any constant $K$, it follows that $\p(Y_{n,n}\geq K \log n)\geq \p(Y_{j,n}\geq K\log n)$ for all $j\in \san$.  We have that 
$$\e(Y_{n,n})=\frac{n}{2}\left(\frac{\frac{2}{\alpha}+\eta}{n}\right)\log n+o(1)=\left(\frac{1}{\alpha}+\frac{\eta}{2}\right)\log n+o(1).$$  An application of Chernoff's bound, see for example \cite[Theorem 2.15]{CL}, gives that for any $\gd>0$, $j\in \san$:
\begin{align*}
\p[Y_{j,n}\geq (1+&\gd)\e(Y_{n,n})]\leq\p\left[Y_{n,n}\geq (1+\gd)\e(Y_{n,n})\right]\\
&=\p\bigg[Y_{n,n}\geq (1+\gd)\bigg(\left[\frac{1}{\alpha}+\frac{\eta}{2}\right]\log n+o(1)\bigg)\bigg]\\
&\leq(1+o(1)) \exp\left\{-\left(\frac{1}{\alpha}+\frac{\eta}{2}\right)(\log n)[(1+\gd)\log(1+\gd)-\gd]\right\}.
\end{align*}
Letting $f(\gd)=\left(\frac{1}{\alpha}+\frac{\eta}{2}\right)[(1+\gd)\log(1+\gd)-\gd]$, we see that $f$ is unbounded and monotonically increasing for $\gd>0$, and so for any $\lambda>0$, an appropriate $\gd_0$ can be chosen  such that $f(\gd_0)=\gl+1$ giving that
$$\p \bigg[Y_{j,n}\geq (1+\gd_0)\bigg(\left[\frac{1}{\alpha}+\frac{\eta}{2}\right]\log n+o(1)\bigg)\bigg]\leq n^{-\gl-1}.$$

Next note that $\rho_{2,n}(j)$ is minimized for $j=\alpha n=(2-\alpha)n$, so that for any constant $K'$, it follows that
$\p(Y_{\alpha n,n}\leq K' \log n)\geq \p(Y_{j,n}\leq K'\log n)$ for all $j\in \san$.  Now $$\e(Y_{\alpha n,n})=\frac{\alpha n}{2}\left(\frac{\frac{2}{\alpha}+\eta}{n}\right)\log n+o(1)=\left(1+\frac{\eta\alpha}{2}\right)\log n+o(1).$$
Another application of Chernoff's bound, see \cite[Theorem 2.17]{CL}, gives then that for any $0\leq \varepsilon\leq e^{-1}$ and $j\in\san$:
\begin{align*}
\p(Y_{j,n}\leq \varepsilon \e(Y_{\alpha n,n}))&\leq\p(Y_{\alpha n,n}\leq \varepsilon \e(Y_{\alpha n,n}))\\
&=\p\bigg(Y_{\alpha n,n}\leq \varepsilon \bigg[\left(1+\frac{\eta\alpha}{2}\right)\log n+o(1)\bigg]\bigg)\\
&\leq(1+o(1))\exp\bigg\{-(1-2\varepsilon+2\varepsilon\log \varepsilon)\left(1+\frac{\eta\alpha}{2}\right)\log n\bigg\}.
\end{align*}
Let $g(\varepsilon)=(1-2\varepsilon+2\varepsilon\log \varepsilon)\left(1+\frac{\eta\alpha}{2}\right)$.  We see that $\lim_{\varepsilon\rightarrow 0}g(\varepsilon)=1+\frac{\eta\alpha}{2}$ so that there exists a $\gamma>0$ and a $\varepsilon_0>0$ such that $g(\varepsilon_0)=1+\gamma$ and 
$$\p(Y_{j,n}\leq \varepsilon_0 \e(Y_{\alpha n,n}))\leq n^{-\gamma-1}$$ for all $j\in\san$.  Next, note that
\begin{align*}
\p(X_n=0)&=\p\left(\cap_{j\in\san} \{Y_{j,n}=\Theta(\log n)\}\right)\\
&=1-\p\left(\cup_{j\in\san} \{Y_{j,n}\neq\Theta(\log n)\}\right)\\
&\geq 1-\sum_{j\in\san} \p(Y_{j,n}\neq\Theta(\log n))\\
&\geq 1-n(n^{-\gamma-1}+n^{-\lambda-1})\\
&=1-n^{-\gamma}-n^{-\lambda}\\
&\rightarrow 1 \text{ as }n\rightarrow\infty,
\end{align*}
which finishes the proof.
\hfill\end{proof}
\noindent{\bf Remark:}
Note that, with the notation as in Theorem 2.1, if for any constants $K,\ \varepsilon>0$ we have  
$$p=p_n:=\sqrt{\frac{K\log^{1+\varepsilon} n}{n}},$$ then 
$$\e(Y_{n,n})=\frac{n}{2}\frac{K\log^{1+\epsilon}n}{n}+o(1)=\Theta(\log^{1+\varepsilon}n).$$  As we have that $Y_{n,n}\approx$Bin$(\frac{n}{2},p^2)$, it follows that 
$$\text{Var}(Y_{n,n})=\Theta(\log^{1+\varepsilon}n),$$  and a simple application of Chebyshev's inequality gives 
$$\p(|Y_{n,n}-\e(Y_{n,n})|\leq \log n)\geq 1-\Theta([\log ^{1-\varepsilon} n]^{-1})\rightarrow 1\text{ as }n\rightarrow\infty.$$  Therefore $\p(X_n=0)\rightarrow 0$ and $\cA_{2,n}$ is not an asymptotic representative 2-basis.

In \cite{AV}, the authors were able to show that if 
$$p=p_n:=\sqrt{\frac{\frac{2}{\alpha}\log n-\frac{2}{\alpha}\log\log n+A_n}{n}}$$ for an arbitrary sequence $A_n=o(\log\log n)$, then 
$$\p(\cA_{2,n}\text{ is an truncated 2-basis})=\begin{cases}
1&\text{ if }A_n\rightarrow\infty\\
0&\text{ if }A_n\rightarrow-\infty\\
\exp\{-2\alpha e^{-\alpha A/2}\}&\text{ if }A_n\rightarrow A.
\end{cases}$$
It follows immediately that if 
$$p=p_n:=\sqrt{\frac{K\log n}{n}},$$ for some $0<K<2$, then with probability converging to 1, $\cA_{2,n}$ will not be a $k$-basis and hence cannot be a representative 2-basis.  At the threshold value 
$$p=p_n\approx\sqrt{\frac{\frac{2}{\alpha}\log n}{n}},$$ we have that the behavior of lower order terms controls whether $\cA_{2,n}$ represents each integer at least once, and so it is reasonable to expect that there are integers that are only represented a few times in the 2-sum set of $\cA_{2,n}$ and therefore $\cA_{2,n}$ will not form a representative 2-basis, though we have no concrete proof of this conjecture beyond our heuristic reasoning.  Note however the similarity between the $p=p_n$ of Theorem \ref{t:2rep} and the $p(z)$ used in \cite{ET} to construct their infinite representative basis.

\section{$k$-Additive Representations}
The problem is complicated further if we consider the representation question in the $k$-additive basis case, as different $k$-sums summing to an integer $j$ are not necessarily disjoint.  We shall begin, as before, by creating the random set $\cA_{k,n}$ by choosing each integer in $[n]\cup\{0\}$ to be in $\cA_{k,n}$ independently with probability $p=p_{n}$.  Fix $\alpha\in(0,1)$, and let $S_k(\alpha,n):=[\alpha n, (k-\alpha)n]$, and for each $j\in S_k(\alpha,n)$, let $Y_{k,n}(j)$ be the number of ways that $j$ can be represented as a $k$-sum of distinct elements of $\cA_{k,n}$.  Let $I_{j,k,n}:=\mathbbm{1}\{Y_{k,n}(j)\neq \Theta(\log n)\}$, so that $X_{k,n}:=\sum_{j\in S_k(\alpha,n)} I_{j,k,n}$ is the number of elements of $S_k(\alpha,n)$ that are not represented order $\log n$ times in the $k$-sum set of $\cA_{k,n}$.  The following theorem will constitute the main result of the section, and the remainder of the section will be dedicated to its proof.
\begt
\label{t:krep}
Let $\varepsilon>0$ be fixed, and let $k\geq 3$ be fixed.  Create the random set $\cA_{k,n}$ by independently picking each integer in $[n]\cup\{0\}$ to be in $\cA_{k,n}$ with probability $$p=p_{n}:=\sqrt[k]{\frac{K\log n}{n^{k-1}}},$$ 
with
$$K=K_{\alpha,k}:=\frac{(4+\varepsilon)(k!)^2}{\alpha^{k-1}}.$$
Then  $$ \lim_{n\rightarrow \infty}P(X_{k,n}=0)=1,$$
so that w.h.p. $\cA_{k,n}$ is an asymptotic representative k-basis as $n\rightarrow\infty$.
\ent
Fix $k\geq3$.  For $1\le l\le k$, define $Y_{l,n}^*(j)$ to be the size of a maximum collection of disjoint representations of $j$ as a $l$-sum of distinct elements of $\mathcal{A}_{k,n}$.  The $Y_{k,n}^*$'s are significantly simpler to work with than the original $Y_{k,n}$'s, as the difficulty presented by overlapping $k$-sums is circumvented.  This idea was exploited to great effect in our motivational paper \cite{ET}.  A few simple calculations yield that for all $i\in [1, (k-\alpha)n]$
$$\e[Y_{l,n}(i)]=O\left(n^{l-1}p^l\right)=O\left(n^{-1+l/k}\right)n^{o(1)}.$$ 
The disjointness lemma (Lemma 1 in \cite{ET}) implies then that for all $l\leq k-1$,
$$\p(Y_{l,n}^*(i)\geq 3k) \leq O\left(n^{-3}\right)n^{o(1)}.$$

We are ready to establish the following lemmata, the first of which is the analogue of Lemma 10 from \cite{ET}:
\begl
\label{lemma:fewsums}
With notation as above, it follows that for all $i\in [1,(k-\alpha)n]$ we have
$$\p\left(Y_{k-1,n}(i)\geq (3k-1)^{k-1}(k-1)!\right)< O\left(n^{-3}\right)n^{o(1)}.$$
\enl
\begin{proof}
We say that $m$ sets form a $\Delta$-system (of size $m$) if they have pairwise the same intersections.   If $Y_{k-1,n}(i)\geq (3k-1)^{k-1}(k-1)!$, then the $\Delta$-system lemma (Lemma 2, \cite{ET}) implies that the set system composed of the $Y_{k-1,n}(i)$ $(k-1)$-sums of $i$ contains a $\Delta$-system of size $3k$, and we shall denote this system via 
$$\{S^{k-1}_1,\ldots,S^{k-1}_{3k}\}$$ with common pairwise intersection set $R$.  Letting $|R|=r\leq k-2$ and letting the sum of elements of $R$ be equal to $m<i$, it follows that if $\widehat{S}^{k-1}_i:=S^{k-1}_i\setminus R$ then 
 $$\{\widehat{S}^{k-1}_1,\ldots,\widehat{S}^{k-1}_{3k}\}$$
 is a system composed of $3k$ disjoint sets of size $k-1-r$ each summing to $i-m$.  The probability of such a system occurring is bounded above by
 $$\p(Y_{k-1-r,n}^*(i-m)\geq 3k) \leq O\left(n^{-3}\right)n^{o(1)}$$
 as desired.
\hfill\end{proof}
The next result is the analogue of Lemma 11 from \cite{ET}:
\begl
\label{lemma:nooverlap}
With notation as above, let $C_k:=(3k-1)^{k-1}{k}!$.  Then for each $j$ in $[\alpha n, (k-\alpha)n]$, we have
$$\p\left(Y_{k,n}(j)\geq C_k Y_{k,n}^*(j)\right)\leq O\left( n^{-2}\right)n^{o(1)}.$$
\enl
\begin{proof}
Slightly abusing notation, we shall write $x\in Y_{k,n}^*(j)$ to mean that $x$ is in one of the maximum collection of disjoint $k$-sums of $j$ counted by $Y_{k,n}^*(j)$.  Then by Lemma 3.2,
\begin{align*}
\p\left(Y_{k,n}(j)\geq C_k Y_{k,n}^*(j)\right)&\leq \p\left(\bigcup_{x\in[0,n]} \{ Y_{k-1,n}(j-x)\ge \frac{C_k}{{k}},\ x\in Y_{k,n}^*(j)\}\right)\\
&\leq \sum_{x\in[0,n]} \p\left(Y_{k-1,n}(j-x)\ge \frac{C_k}{{k}},\ x\in Y_{k,n}^*(j)\right)\\
&\leq \sum_{x\in[0,n]} \p\left(Y_{k-1,n}(j-x)\ge  \frac{C_k}{{k}}\right)\\
&\leq \sum_{x\in[0,n]} O\left( n^{-3}\right)n^{o(1)}\\
&=O\left( n^{-2}\right)n^{o(1)},
\end{align*}
as desired.
\hfill\end{proof}
Next, we shall use Talagrand's inequality (see Section 7.7, \cite{as}) to show that $Y_{k,n}^*(j)=\Theta(\log n)$ with high probability for all $j\in[\alpha n,(k-\alpha)n]$.  Towards that end, we prove
\begl
\label{lemma:median}
For some constant $B_{k}\in[-40\sqrt{k},40\sqrt{k}]$ we have that
$$\text{\bf Med}(Y_{k,n}^*(j))=\e(Y_{k,n}^*(j))+B_{k}\sqrt{\e(Y_{k,n}^*(j))}.$$
\enl
\begin{proof}
First note that $Y_{k,n}^*(j)$ can be written as a function 
$$Y_{k,n}^*(j)=f(J_0,J_1,\ldots,J_n)$$ of the indicator variables
$$J_i:=\begin{cases}
1&\text{ if }i\in \cA\\
0&\text{ else.}
\end{cases}$$
As the $k$-sums counted by $Y_{k,n}^*(j)$ are disjoint, the function $f(\cdot)$ is  one-Lipschitz.  Also note that $f$ is $h-$certifiable with $h(s)=ks$, since if $f(J_0,\ldots,J_n)\geq s$, there exists $s$ disjoint $k$-sums of $j$ present in $\cA$, and any other realization of $\cA$ with those $sk$ $J_i$'s equal to 1 has $f\geq s$ as well.  It immediately follows from Fact 10.1 in \cite{Mol} that
$$\left|\e(Y_{k,n}^*(j))-\text{\bf Med}(Y_{k,n}^*(j))\right|\leq 40\sqrt{k\e(Y_{k,n}^*(j))}.$$  This completes the proof.
\hfill\end{proof}
\noindent Next we prove
\begl
\label{lem:expbnd}
With $p=p_n$ defined as in Theorem 3.1 and notation as above,
$$\e\left(Y^*_{k,n}(j)\right)\leq \e\left(Y_{k,n}(j)\right)\leq\e\left(Y^*_{k,n}(j)\right)+o(1).$$
\enl
\begin{proof}
Let $W_{j,k,n}$ be the number of overlapping pairs of $k$-sums in the set of all $k$-sums of $j$ using elements of $\cA$.  Then, as each $k$-sum not in $Y^*_{k,n}(j)$ must intersect with at least one of the $k$-sums of $Y^*_{k,n}(j)$, we have that
$$Y^*_{k,n}(j)\leq Y_{k,n}(j)\leq Y^*_{k,n}(j)+W_{j,k,n}.$$  Note that (writing $\sum_{l,*}$ to be the sum over all overlapping pairs of $k$-sums of $j$ using elements of $\cA$ with overlap of size $l$)
\begin{align*}
\e(W_{j,k,n})&=\sum_{l=1}^{k-1}\sum_{l,*} p^{2k-l}\\
&=\sum_{l=1}^{k-1}O\left(n^{2k-l-2}p^{2k-l}\right)\\
&=\sum_{l=1}^{k-1}O\left(n^{-l/k}(\log n)^{(2k-l)/k}\right)\\
&=O\left(n^{-1/k}(\log n)^{(k+1)/k}\right)\\
&=o(1),
\end{align*}
as desired.
\hfill\end{proof}

\begt
\label{thm:ub}
With $p$ defined as in Theorem 3.1, there exists constants $\gamma_j>0$ and $\xi>0$ such that 
$$\p\left(Y_{k,n}^*(j)\leq \gamma_j \log(n)+O(\sqrt{\log n})\right)\leq 2n^{-1-\xi}.$$
\ent
\begin{proof}
In \cite{AV} it was shown that the number $\rho_{k,n}(j)$ of (not necessarily disjoint) $k$-sums of distinct elements of $[n]\cup\{0\}$ summing to $j$, for $j\in[\alpha n,(k-\alpha)n]$, is bounded below by
$$\rho_{k,n}(j)\geq (1+o(1))\frac{(\alpha n)^{k-1}}{k!(k-1)!}.$$ 
It is immediate that $\rho_{k,n}(j)=O(n^{k-1})$, so that there exists a constant $C(j)\geq \frac{\alpha^{k-1}}{k!(k-1)!}$ such that $\rho_{k,n}(j)=C(j)(1+o(1))n^{k-1}$.  From Lemma 3.5, we have then that for $j\in[\alpha n,(k-\alpha)n]$,
\begin{align*}
\e[Y^*_{k,n}(j)]&= (1+o(1))C(j)n^{k-1}p^k+o(1)\\
&=(1+o(1))C(j) K_{\alpha,k}\log n+o(1),
\end{align*}
and so Lemma 3.4 gives us that
\begin{equation}
\label{med}
\text{\bf Med}(Y_{k,n}^*(j))=(1+o(1))C(j) K_{\alpha,k}\log n+O\left(\sqrt{\log n}\right).
\end{equation}
Talagrand's inequality (see Theorem 7.7.1 in \cite{as}) gives us that for all $t,m>0$ (where $h(s)=ks$ is the aforementioned certification function for $Y_{k,n}^*(j)=f(J_0,J_1,\ldots,J_n)$):
$$\p\left(Y_{k,n}^*(j)\leq m-t\sqrt{h(m)}\right)\p\left(Y_{k,n}^*(j)\geq m\right)\leq e^{-t^2/4}.$$  Let $t=\sqrt{(4+4\xi)\log n}$, and 
$m=\text{\bf Med}(Y_{k,n}^*(j))$ to see that
$$\p\left(Y_{k,n}^*(j)\leq\text{\bf Med}(Y_{k,n}^*(j))-\sqrt{(4+4\xi)\log n}\sqrt{k \text{\bf Med}(\rho_{k,n}^*(j))}\right)\leq 2 n^{-1-\xi}.$$
Using (\ref{med}), we see then that
\begin{align*}
\p\bigg(Y_{k,n}^*(j)\leq (1+o(1))\bigg[C(j)K_{\alpha,k}&-\sqrt{4+4\xi}\sqrt{kC(j)K_{\alpha,k}} \bigg]\log n+\cdots\\
&\cdots+O(\sqrt{\log n})\bigg)\leq 2n^{-1-\xi}.
\end{align*}
Letting $\gamma_j=C(j)K_{\alpha,k}-\sqrt{4+4\xi}\sqrt{kC(j)K_{\alpha,k}}$, then for any $\xi<\varepsilon/4$ (where this is the $\varepsilon$ from the definition of $K_{\alpha,k}$) it follows from
$C(j)\geq \frac{\alpha^{k-1}}{k!(k-1)!}$ that $\gamma_j>0$ as desired.
\hfill\end{proof}
\begt
\label{thm:lb}
With $p=p_n$ defined as in Theorem 3.1, for each $j$ there exists a constant $\delta_j>0$ such that
$$\p\left(Y_{k,n}^*(j)\geq \delta_j \log(n)+O(\sqrt{\log n})\right)\leq  2 n^{-5/4}.$$
\ent
\begin{proof}
We will again use Talagrand's inequality, but we shall now set 
$$m-t\sqrt{km}=\text{\bf Med}(Y^*_{k,n}(j)).$$
Solving for $m$, we get that
$$m=\left(\frac{t\sqrt k}{2}+\frac{1}{2}\sqrt{kt^2+4\text{\bf Med}((Y^*_{k,n}(j))}\right)^2.$$
As in the proof of Theorem 3.6, we have that 
$$\text{\bf Med}[(Y^*_{k,n}(j))=(1+o(1))C(j)K_{\alpha,k}\log n+O\left(\sqrt{\log n}\right),$$
with $C(j)\geq \frac{\alpha^{k-1}}{k!(k-1)!}$ so that
$$m=\left(\frac{t\sqrt k}{2}+\frac{1}{2}\sqrt{kt^2+4(1+o(1))C(j)K_{\alpha,k}\log n+O\big(\sqrt{\log n}\big)}\right)^2.$$
Let $t=\sqrt{ 5\log n}$ to arrive at
$m=\delta_j \log n+O(\sqrt{\log n})$
for some constant $\delta_j$.  Apply Talagrand's inequality to see that
$$\p\left(Y_{k,n}^*(j)\geq\delta_j \log n+O\big(\sqrt{ \log n}\big)\right)\leq 2 n^{-5/4}$$ as desired.
\end{proof}
We are now ready to prove our main result:
\begin{proof}[Proof of Theorem 3.1:]
Let $$\gamma_n:=\min_{j\in S_k(\alpha,n)}\gamma_j, \text{ and }\delta_n:=\max_{j \in S_k(\alpha,n)}\delta_j.$$  Note that there exist strictly positive finite functions  $g_1(k)$, $g_2(k)$, $g_3(k)$ and $g_4(k)$ of $k$, such that for all $n$, $ g_1(k)<\gamma_n<g_2(k)$ and $g_3(k)< \delta_n< g_4(k)$.  It follows that as $n\rightarrow\infty$, we have 
$$0<\lim_{n\rightarrow\infty} \gamma_n<\infty,\text{ and } 0<\lim_{n\rightarrow\infty}\delta_n<\infty.$$  
It follows from Theorems 3.6 and 3.7 that there exists a $\xi>0$ such that for all $j\in S_k(\alpha,n)$,
$$\p\left(Y_{k,n}^*(j)\leq\gamma_n \log n+O\big(\sqrt{\log n}\big)\right)\leq 2n^{-1-\xi},$$
$$\p\left(Y_{k,n}^*(j)\geq\delta_n \log n+O\big(\sqrt{\log n}\big)\right)\leq 2n^{-5/4}.$$
It follows immediately that for any constant $c$,
$$\p\left(Y_{k,n}^*(j)\leq c\right)\geq \p\left(Y_{k,n}(j)\leq c\right),$$ and hence there exists a $\xi>0$ such that for all $j\in S_k(\alpha,n)$,
$$\p\left(Y_{k,n}(j)\leq\gamma_n \log n+O\big(\sqrt{\log n}\big)\right)\leq 2n^{-1-\xi}.$$
Next note that by Lemma 3.3,
\begin{align*}
&\p\left(Y_{k,n}(j)\geq  C_k\delta_n\log n+O(\sqrt{\log n}\right)\\
&=\p\left(Y_{k,n}(j)\geq  C_k\delta_n\log n+O\big(\sqrt{\log n}\big),\  Y_{k,n}(j)\geq C_k Y^*_{k,n}(j)\right)\\
&\hspace{5mm}+\p\left(Y_{k,n}(j)\geq  C_k\delta_n\log n+O\big(\sqrt{\log n}\big),\  Y_{k,n}(j)<C_k Y^*_{k,n}(j)\right)\\
&\leq O(n^{-2})n^{o(1)}+\p\left(Y^*_{k,n}(j)\geq \delta_n\log n+O\big(\sqrt{\log n}\big)\right)\\
&=O(n^{-2})n^{o(1)}+2n^{-5/4}.
\end{align*}
Therefore, defining the event
$$A_j:= \{Y_{k,n}(j)\geq C_k \delta_n\log n+O\big(\sqrt{\log n}\big)\}\cup\{ Y_{k,n}(j)\leq \gamma_n\log n+O\big(\sqrt{\log n}\big)\}$$
\begin{align*}
\p(X_{k,n}\geq 1)&\leq \p\left(\cup_{j} A_j\right)\\
&\leq \sum_j \p\left(Y_{k,n}(j)\geq C_k \delta_n\log n+O\big(\sqrt{\log n}\big)\right)\\
&\hspace{15mm}+\sum_j \p\left(Y_{k,n}(j)\leq \gamma_n\log n+O\big(\sqrt{\log n}\big)\right)\\
&\leq kn \left( O(n^{-2})n^{o(1)}+2n^{-5/4}+2n^{-1-\xi}\right)\\
&=O(n^{-\xi})=o(1),
\end{align*}
and $\p(X_{k,n}=0)\rightarrow 1$ as $n\rightarrow\infty$ as desired.
\end{proof}
\noindent{\bf Remarks:}  (i) If we consider representations of integers in $[\alpha n,(k-\alpha)n]$ using $k$ integers from $\cA_{k,n}$  that are not necessarily distinct, we can prove a result similar to Theorem 3.1.  We skip the details.

\noindent(ii) In \cite{AV}, the authors showed that if 
$$p:=\sqrt[k]{\frac{\frac{k!(k-1)!}{\alpha^{k-1}}\log n-\frac{k!(k-1)!}{\alpha^{k-1}}\log\log n+A_n}{n^{k-1}}}$$ for $A_n=o(\log\log n)$,
then 
$$\p(\cA_{k,n}\text{ is an asymptotic k-basis})\rightarrow\begin{cases}
1&\text{ if }A_n\rightarrow\infty\\
0&\text{ if }A_n\rightarrow-\infty\\
\exp\left\{-\frac{2\alpha}{k-1}e^{\frac{-A\alpha^{k-1}}{(k!(k-1)!)}}\right\}&\text{ if }A_n\rightarrow A<\infty.
\end{cases}$$
Therefore if we choose elements to be in $\cA_{k,n}$ with probability
$$p=\sqrt[k]{\frac{C\log n}{n^{k-1}}}$$
for $C<\frac{k!(k-1)!}{\alpha^{k-1}}$,
then w.h.p. $\cA_{k,n}$ is not a $k$-basis, and so w.h.p. it is not a representative $k$-basis.  When 
$C=\frac{k!(k-1)!}{\alpha^{k-1}}$, the behavior of $\cA_{k,n}$ as a $k$-basis hinges on the behavior of lower order terms, and so again it is reasonable to expect some integers to be represented only a few times as $k$-sums of elements of $\cA_{k,n}$.  We would expect then that $\cA_{k,n}$ is not a representative $k$-basis w.h.p.  As the constant $C$ increases to $K_{\alpha,k}$, our random set becomes a representative $k$-basis w.h.p. as $n\rightarrow\infty$.  We haven't yet established any threshold behavior when $\frac{k!(k-1)!}{\alpha^{k-1}}\leq C<K_{\alpha,k}$, leaving the door open for future research.

\section{Further Research}  It would be interesting to work out the asymptotics for when $\cA$ becomes a {truncated $\phi(n)$-representative $k$-basis} for $\phi(n)=o(n^\varepsilon)$ for $\phi(n)$ other than $\log n$, and to what extent the linearity of the target sumset can be relaxed from $[\alpha n,(k-\alpha)n]$ to perhaps $[\theta(n), kn-\theta(n)]$.  We expect the most difficult challenges to present themselves for $\phi(n)=o(\log n)$.  In a similar vein, improvements in Theorems 2.1 and 3.1 would be most instructive.  In particular, how close to the additive basis constant $\frac{k!(k-1)!}{\alpha^{k-1}}$ can we force the constant in Theorem 3.1, as discussed in greater detail in the previous paragraph?

\section{Acknowledgments} The research of all four authors was supported by NSF Grant 1004624.  The research of VL was supported by the Acheson~J.~Duncan Fund for the Advancement of 
Research in Statistics, and by U.S. Department of Education GAANN grant P200A090128.


\begin{thebibliography}{99}
\bibitem{as} N.~Alon and J.~Spencer (1992).  {\it The Probabilistic Method.} Wiley, New York.
\bibitem{bhj} A.~Barbour, L.~Holst, and S.~Janson (1992).  {\it Poisson Approximation.}  Oxford University Press.
\bibitem{CL} F.~Chung and L.~Lu (2006).  {\it Complex Graphs and Networks.}  NSF-CBMS Lecture Notes Series, American Mathematical Society, Providence.
\bibitem{erdos} P.~Erd\H os (1956). Problems and results in additive number theory, in {\it Colloque sur la Th\'eorie des Nombres} (CBRM), Bruxelles, 127--137.
\bibitem{ET} P.~Erd\H os and P.~Tetali (1990). Representations of integers as the sum of $k$ terms, {\it Rand. Structures Algorithms} {\bf 1}, 245--261.
\bibitem{gjlr} A.~Godbole, S.~Janson, N.~Locantore, and R.~Rapoport (1999).  Random Sidon sequences,  {\it J. Number Theory} {\bf 75}, 7--22.
\bibitem{AV} A.~Godbole, C-M Lim, V.~Lyzinski, and N.Triantafillou (2013).  Sharp threshold asymptotics for the emergence of additive bases, {\it Integers:  Electronic Journal of Combinatorial Number Theory} {\bf 13}, Paper \#
\bibitem{gn} C.~G\"unt\"urk and M.~Nathanson (2006). A new upper bound for finite additive bases, {\it Acta Arith.}
{\bf 124}, 235--255.
\bibitem{hh} N.~H\"ammerer and G.~Hofmeister (1976).  Zu einer Vermutung von Rohrbach, {\it J. Reine Angew. Math.}
{\bf 286-287}, 239--247.
\bibitem{janson} S.~Janson (1990). Poisson approximation for large deviations, {\it Rand. Structures Algorithms} {\bf 1}, 221--229.
\bibitem{Mol} M.~Molloy and B.~Reed (2001).  {\it Graph Coloring and the Probabilistic Method.}  Springer Verlag, New York.
\bibitem{mrose} A. Mrose (1979).  Untere Schranken f\"ur die Reichweiten von Extremalbasen fester Ordnung, {\it Abh.
Math. Sem. Univ. Hamburg} {\bf 48}, 118--124.
\bibitem{sandor} C.~Sandor (2007).  Random $B_h$ sets and additive bases in ${\mathbb Z}_n$, {\it Integers: Electronic Journal of Combinatorial Number Theory} {\bf 7}, Paper \#A32.
\bibitem{yu} G.~Yu (2009).  Upper bounds for finitely additive 2-bases, {\it Proc. Amer. Math. Soc.} {\bf 137}, 11--18.
\bibitem{ya} A.~Yadin (2009).  When Do Random Subsets Decompose a Finite Group? {\it Israel J.~Math.} {\bf 174}, 203--219.
\end{thebibliography}
\end{document}